\documentclass{article}

\newtheorem{thm}{Theorem}
\newtheorem{defn}[thm]{Definition}
\newtheorem{lem}[thm]{Lemma}
\newtheorem{cor}[thm]{Corollary}

\newenvironment{pfmetric}{\medskip \noindent
{\bf Proof of Theorem \ref{thm;metric}.}}{\hfill \rule{.5em}{1em}\\}
\newenvironment{pfinequality}{\medskip \noindent
{\bf Proof of Lemma \ref{lem;inequality}.}}{\hfill \rule{.5em}{1em}\\}
\newenvironment{pfeuclid}{\medskip \noindent
{\bf Proof of Corollary \ref{cor;euclid}.}}{\hfill \rule{.5em}{1em}\\}

\newcommand{\mb}{\mathbf}

\newcommand{\ep}{\epsilon}

\title{Locally Euclidean metrics on $\mb{R}^3$}
\author{Young Deuk Kim \\ School of Mathematical Sciences\\
Seoul National University\\ Seoul, 151-747, Korea\\(ydkimus@yahoo.com)}
\date{\today}

\begin{document}
\maketitle

\begin{abstract}
For all $0<t\leq 1$, we define a locally Euclidean metric $\rho_t$ on
$\mb{R}^3$. These metrics are invariant under Euclidean isometries and, if
$t$ increases to $1$, converge to the Euclidean metric $d_E$.
This research is motivated by expanding universe.\\

\noindent
{\footnotesize{\textbf{key words.}}} locally Euclidean metric\\

\noindent
{\footnotesize{\textbf{PACS number(s).}}} 98.80.Jk \\
{\footnotesize{\textbf{Mathematics Subject Classifications (2000).}}} 85A40,
57M50
\end{abstract}

\section{The metric $\rho_t$}

Let $d_E$ denote the Euclidean metric on $\mb{R}^3$. As in \cite{KIM},
a nonnegative function $d:\mb{R}^3\times\mb{R}^3\to\mb{R}$ is called a
{\em locally Euclidean metric} if
\begin{enumerate}
\item[(i)]
$d(P,Q)=0$ if and only if $P=Q$
\item[(ii)]
$d(P,Q)=d(Q,P)$ for all $P,Q\in \mb{R}^3$
\item[(iii)]
$d(P,Q)+d(Q,R)\geq d(P,R)$ for all $P,Q,R\in \mb{R}^3$
\item[(iv)]
For all $P\in\mb{R}^3$, there exists $r>0$ such that
$$d(Q,R)=d_E(Q,R)\quad\mbox{for all }
Q,R\in B_r(P)=\{S\in\mb{R}^3\mid d(P,S)<r\}.$$
\end{enumerate}

\noindent
Let $S^2\subset\mb{R}^3$ be the unit 2-sphere with center $O=(0,0,0)$ and
radius $1$. A locally Euclidean metric on $S^2$ is defined in the same way
as on $\mb{R}^3$. In the following theorem of the author, $-P$ is the antipodal
point of $P\in S^2$ and $0\leq\alpha<\pi/4$.
\begin{thm}[\cite{KIM}]\label{thm;kim}
Suppose that $0<t\leq 1$ and
$\alpha=\sin^{-1}\left({\sqrt{2-t^2}-t}\over 2\right)$.
For all $P,Q\in S^2$ we define
\begin{eqnarray*}
d_t(P,Q)=\left\{
\begin{array}{lll}
d_E(P,Q) &\mbox{ if }& \angle POQ\leq \pi-2\alpha\\
2t+d_E(-P,Q) &\mbox{ if }& \angle POQ>\pi-2\alpha.
\end{array}
\right.
\end{eqnarray*}
Then $d_t$ is a locally Euclidean metric on $S^2$
which is invariant under any Euclidean isometry and
\begin{equation}\label{eq;t1}
td_E(P,Q)\leq d_t(P,Q)\leq d_E(P,Q).
\end{equation}
\end{thm}

\indent
Notice that
\begin{equation}
\mbox{if}\quad d_E(P,Q)=d_E(R,S)\quad \mbox{then}\quad
d_t(P,Q)=d_t(R,S)\label{eq;dteuclidean}
\end{equation}
and
\begin{equation}
\mbox{if}\quad d_t(P,Q)<2t\quad \mbox{then}\quad
d_t(P,Q)=d_E(P,Q).\label{eq;te}
\end{equation}

\indent
We will make use of the metric $d_t$ to define a locally Euclidean metric
$\rho_t$ on $\mb{R}^3$. Suppose that $P,Q\in\mb{R}^3$ and $n\in\mb{N}$. Let
\begin{eqnarray*}
\Gamma^n_{P,Q}&=&\{\, (X_0,X_1,\cdots,X_n)\mid X_0=P,\ X_n=Q\\
&&\mbox{and}\quad X_i\in\mb{R}^3,\ 0\leq d_E(X_{i-1},X_i)\leq 2\
\mbox{ for all }1\leq i\leq n\, \}
\end{eqnarray*}
and
$$\Gamma_{P,Q}=\bigcup_{n\in\mb{N}}\Gamma^n_{P,Q}.$$
Notice that $\Gamma_{P,Q}\neq\emptyset$ for all $P,Q\in\mb{R}^3$ and
if $(X_0,X_1,\cdots,X_n)\in\Gamma_{P,Q}$ then
$(X_n,X_{n-1},\cdots,X_0)\in\Gamma_{Q,P}$.

\indent
Suppose that $(X_0,X_1,\cdots,X_n)\in\Gamma_{P,Q}$ and $1\leq i\leq n$.
There exists an isometric embedding $\ep_i:(S^2,d_E)\to(\mb{R}^3,d_E)$
such that $X_{i-1},X_i\in\ep(S^2)$. Notice that
\begin{equation}\label{eq;epi}
d_t\left({\ep_i}^{-1}(X_{i-1}),{\ep_i}^{-1}(X_i)\right)
\end{equation}
does not depend on the choice of $\ep_i$ because $d_t$ is invariant under
Euclidean isometries. Therefore for the sake of simplicity, we will
simply write $$d_t(X_{i-1},X_i)$$ to denote eq. (\ref{eq;epi}).
In the following definition, the infimum is taken for all elements of
$\Gamma_{P,Q}$.
\begin{defn}\label{def;metric}
For $0<t\leq 1$ and $P,Q\in\mb{R}^3$, we define
$$\rho_t(P,Q)=\inf_{(X_0,X_1,\cdots,X_n)\in\Gamma_{P,Q}}
\sum^n_{i=1}d_t(X_{i-1},X_i).$$
\end{defn}

\indent
We will prove the following lemma and theorem in the next section.
\begin{lem}\label{lem;inequality}
$td_E(P,Q)\leq \rho_t(P,Q)\leq d_E(P,Q)$ for all $P,Q\in\mb{R}^3$.
\end{lem}

\begin{thm}\label{thm;metric}
$\rho_t$ defines a locally Euclidean metric on $\mb{R}^3$.
\end{thm}

\indent
Suppose that $\phi:\mb{R}^3\to\mb{R}^3$ is an Euclidean isometry.
Notice that
$$(X_0,X_1,\cdots,X_n)\in\Gamma_{P,Q}\mbox{ if and only if }
(\phi(X_0),\phi(X_1),\cdots,\phi(X_n))\in\Gamma_{\phi(P),\phi(Q)}.$$
Therefore, from eq. (\ref{eq;dteuclidean}) we can show that
$$\rho_t(P,Q)=\rho_t(\phi(P),\phi(Q)).$$

\indent
By Lemma \ref{lem;inequality} if $t$ increases to $1$, then
$$\rho_t(P,Q)\ \mbox{converges to}\ d_E(P,Q)\quad
\mbox{for all } P,Q\in\mb{R}^3.$$

\indent
The following corollary which will be proven in the next section implies
that $\rho_t$ and $d_E$ are same locally.
\begin{cor}\label{cor;euclid}
If $d_E(P,Q)<2t$ then $\rho_t(P,Q)=d_E(P,Q)$ for all $P,Q\in\mb{R}^3$.
\end{cor}

\indent
The construction of $\rho_t$ is motivated by {\em expanding universe}.
See \cite{BARROW}, \cite{ELLIS} and \cite{KAFATOS} as references on this topic.
Notice that $(\mb{R}^3,\rho_t)$ is homogeneous, isotropic and homeomorphic to
the Euclidean space $(\mb{R}^3,d_E)$ for all $0< t\leq 1$.
From the proof of Theorem \ref{thm;metric} in the next section, notice also 
that if $t$ increases, so is the Euclidean region $B_t(P)$ for all 
$P\in\mb{R}^3$.

\indent The Einstein's field equations are expressed in terms of
Riemannian metrics (see \cite{NORTON}, \cite{OHANIAN} and \cite{RAY}
for details). If a Riemannian metric space is locally Euclidean,
then it is the Euclidean space itself. Therefore the construction of
this paper is impossible with Riemannian metrics. In a preprint of
the author(\cite{KIM-preprint}), there exists a similar construction of 
pseudo metrics on the closed ball $B^3$.

\section{Proof of the lemma, theorem and corollary}

\begin{pfinequality}
Suppose that $P,Q\in\mb{R}^3$. Choose $m\in\{0\}\cup\mb{N}$ such that
$m\leq d_E(P,Q)<m+1$.
Find $(X_0,X_1,\cdots,X_{m+1})\in\Gamma_{P,Q}$ such that
$$X_i\in\overline{PQ}\quad\mbox{and}\quad
d_E(X_{i-1},X_i)=\frac{1}{m+1}d_E(P,Q)\quad\mbox{for all }1\leq i\leq m+1,$$
where $\overline{PQ}$ is the Euclidean segment from $P$ to $Q$.
From eq (\ref{eq;t1}), we have
\begin{eqnarray*}
&&\rho_t(P,Q)\leq\sum_{i=1}^{m+1}d_t(X_{i-1},X_i)
\leq\sum_{i=1}^{m+1}d_E(X_{i-1},X_i)\\
&&\qquad\qquad\qquad\qquad\qquad\quad
 =\sum_{i=1}^{m+1}\frac{1}{m+1}d_E(P,Q)=d_E(P,Q).
\end{eqnarray*}
\indent
For all $(X_0,X_1,\cdots,X_n)\in\Gamma_{P,Q}$, from eq. (\ref{eq;t1}) we have
$$\sum_{i=1}^{n}d_t(X_{i-1},X_i)\geq
\sum_{i=1}^{n}td_E(X_{i-1},X_i)=
t\sum_{i=1}^{n}d_E(X_{i-1},X_i)\geq
td_E(P,Q).$$
Thus $\rho_t(P,Q)\geq td_E(P,Q)$.
\end{pfinequality}

\begin{pfmetric}
Since $d_t$ is nonnegative, so is $\rho_t$. By Lemma \ref{lem;inequality},
it is trivial to show that $\rho_t(P,Q)=0$ if and only if $P=Q$.
Since $(X_0,X_1,\cdots,X_n)\in\Gamma_{P,Q}$ if and only if
$(X_n,X_{n-1},\cdots,X_0)\in\Gamma_{Q,P}$, from eq. (\ref{eq;dteuclidean})
we have 
$$\rho_t(P,Q)=\rho_t(Q,P)\quad\mbox{for all }P,Q\in\mb{R}^3.$$

\indent
For all $(X_0,X_1,\cdots,X_n)\in\Gamma_{P,Q}$ and
$(Y_0,Y_1,\cdots,Y_m)\in\Gamma_{Q,R}$, we have
$$(X_0,X_1,\cdots,X_n=Y_0,Y_1,\cdots,Y_m)\in\Gamma_{P,R}.$$
Therefore
$$\rho_t(P,R)\leq \sum_{i=1}^{n}d_t(X_{i-1},X_i)+
\sum_{i=1}^{m}d_t(Y_{i-1},Y_i)$$
and hence
\begin{equation}\label{eq;triangle}
\rho_t(P,R)\leq\rho_t(P,Q)+\rho_t(Q,R)\quad\mbox{for all }P,Q,R\in\mb{R}^3.
\end{equation}

\indent
Suppose that $P\in\mb{R}^3$ and
$Q,R\in B_t(P)=\{S\in\mb{R}^3\mid \rho_t(P,S)<t\}$.
From eq. (\ref{eq;triangle}), we have
$\rho_t(Q,R)\leq \rho_t(Q,P)+\rho_t(P,R)<2t$. Therefore there exists
$(X_0,X_1,\cdots,X_n)\in\Gamma_{Q,R}$ such that
$$\sum_{i=1}^{n}d_t(X_{i-1},X_i)<2t.$$
From eq. (\ref{eq;te}), we have
$d_t(X_{i-1},X_i)=d_E(X_{i-1},X_i)$ for all $1\leq i\leq n$.
Therefore
\begin{eqnarray*}
&&\rho_t(Q,R)=\inf_{\Gamma_{Q,R}}
\sum^n_{i=1}d_t(X_{i-1},X_i)
=\inf_{\Gamma_{Q,R}(t)}\sum^n_{i=1}d_t(X_{i-1},X_i)\\
&&\qquad\qquad\qquad\qquad\qquad\qquad
=\inf_{\Gamma_{Q,R}(t)}\sum^n_{i=1}d_E(X_{i-1},X_i)\geq d_E(Q,R),
\end{eqnarray*}
where
$$\Gamma_{Q,R}(t)=\{\, (X_0,X_1,\cdots,X_n)\in\Gamma_{P,Q}\, \mid \,
\sum_{i=1}^{n}d_t(X_{i-1},X_i)<2t\, \}.$$
Thus from Lemma \ref{lem;inequality}, we have $\rho_t(Q,R)=d_E(Q,R)$.
\end{pfmetric}

\begin{pfeuclid}
Suppose that $d_E(P,Q)<2t$.
From Lemma \ref{lem;inequality}, we have
\begin{equation}\label{eq;2t}
\rho_t(P,Q)<2t.
\end{equation}
In the last paragraph of the proof of Theorem \ref{thm;metric}, we showed
that eq. (\ref{eq;2t}) implies $\rho_t(P,Q)=d_E(P,Q)$.
\end{pfeuclid}



\begin{thebibliography}{99}

\bibitem{BARROW}
J.D. Barrow, \textit{Relativistic cosmology}, The early universe
(eds. W. G. Unruh and G. W. Semenoff), 125-201, Reidel, Dordrecht, 1988.

\bibitem{ELLIS}
G.F.R. Ellis, \textit{The expanding universe: A history of cosmology 
from 1917 to 1960}, Einstein and the history of general relativity
(eds. D. Howard and J. Stachel), 367-431, Birkh\"auser, Boston, 1989.

\bibitem{KAFATOS}
M. Kafatos, \textit{The problem of observation in cosmology and the big bang},
Gravitation and cosmology: From the Hubble radius to the Plank scale
(eds. R.L. Amoroso, G. Hunter, M. Kafatos and J.P. Vigier), 65-80,
Kluwer Academic Publishers, Dordrecht, 2002.

\bibitem{KIM-preprint}
Y.D. Kim, \textit{A family of pseudo metrics on $B^3$ and its application},
preprint, arXiv:math.GN/0201077.

\bibitem{KIM}
Y.D. Kim, \textit{Locally Euclidean metrics on $S^2$ in which some open balls
are not connected}, Rocky Mountain Journal of Mathematics {\bf 36}(2006), no.6,
1927-1935.

\bibitem{NORTON}
J. Norton, \textit{How Einstein found his field equations, 1912-1915},
Einstein and the history of general relativity (eds. D. Howard and J. Stachel),
101-159, Birkh\"auser, Boston, 1989.

\bibitem{OHANIAN}
H.C. Ohanian, Gravitation and spacetime, W. W. Norton \& Company,
New York, 1976.

\bibitem{RAY}
C. Ray, \textit{The cosmological constant: Einstein's greatest mistake?},
Stud. Hist. Philos. Sci. {\bf 21}(1990), no.4, 589-604.

\end{thebibliography}
\end{document}